\documentclass[12pt]{amsart}
\usepackage{amsmath,amssymb,amscd,amsxtra,amsthm}
\usepackage{latexsym}
\usepackage[dvips]{graphics,epsfig}

\newcommand{\SG}{\text{SG}}

\newcommand{\R}{\mathbb{R}}

\newcommand{\set}[1]{\{#1\}}

\newcommand{\nm}[1]{\|{#1}\|}

\newcommand{\biggparen}[1]{\biggl(#1\biggr)}
\newcommand{\ip}[2]{\langle#1,#2\rangle}

\newcommand{\en}[1]{\mathcal{E}(#1,#1)}

\newcommand{\vstr}[1][3]{\rule{0ex}{#1ex}}
\newcommand{\negsp}[1][20]{\mspace{-#1mu}}
\newcommand{\evald}[2][]{\ensuremath{\negsp[4]\left.\vstr[2.0] \right|_{#2}^{#1}}} 

\DeclareMathOperator{\trace}{Trace}%

\newtheorem{thm}{Theorem}
\newtheorem{lemma}{Lemma}
\newtheorem{prop}{Proposition}

\theoremstyle{remark}
\newtheorem{rem}{Remark}
\theoremstyle{definition}
\newtheorem{deft}{Definition}

\begin{document}
\title{Szeg\"o limit theorems on the Sierpi\'nski gasket}

\author{Kasso A.~Okoudjou}
\address{Kasso A.~Okoudjou\\
Department of Mathematics\\
University of Maryland\\
College Park, MD 20742-4015, USA} \email{kasso@math.umd.edu}

\author{Luke G.~Rogers}
\address{Luke G.~Rogers\\
Department of Mathematics\\
University of Connecticut\\
Storrs, CT 06269-3009, USA} \email{rogers@math.uconn.edu}

\author{Robert S.~Strichartz} 
\thanks{Research of the third author supported in part by the National Science Foundation, grant DMS-065440.}

\address{Robert S.~Strichartz\\
Department of Mathematics\\
Malott Hall\\
Cornell University, Ithaca, NY 14853-4201, USA} \email{str@math.cornell.edu}


\date{\today}


\begin{abstract}
We use the existence of localized eigenfunctions of the Laplacian on the Sierpi\'nski gasket
($\SG$) to formulate and prove  analogues of  the strong Szeg\"o limit theorem in this fractal
setting. Furthermore, we recast some of our results in terms of equally distributed sequences.
\end{abstract}

\maketitle \pagestyle{myheadings} \thispagestyle{plain} \markboth{K. A. OKOUDJOU, L. G. ROGERS AND
R. S. STRICHARTZ}{SZEG\"O LIMIT THEOREM}

\section{Introduction}\label{sec1}
Let $P_n$ be the orthogonal projection of $L^{2}([0, 2\pi)$ onto the linear subspace spanned by the
functions $\{e^{i m \theta}: 0\leq m\leq n; 0\leq \theta < 2\pi\}$. For any function $f$ defined on
$[0, 2\pi)$, let $[f]$ be the linear operator corresponding to multiplication by $f$.

In $1952$, G.~Szeg\"o proved that for a positive function $f \in \mathcal{C}^{1+\alpha}$ where
$\alpha > 0$, the following holds
\begin{equation}\label{oriszeg}
\lim_{n \to \infty}\tfrac{1}{n+1}\log \det P_{n}[f]P_{n}=\tfrac{1}{2\pi}\int_{0}^{2\pi}\log f
(\theta)\, d\theta.
\end{equation} Equivalently,~\eqref{oriszeg} can  be expressed as
\begin{equation}\label{oriszegbis}
\lim_{n \to \infty}\tfrac{1}{n+1}\trace \log P_{n}[f]P_{n}=\tfrac{1}{2\pi}\int_{0}^{2\pi}\log f
(\theta)\, d\theta.
\end{equation}
The above result is known as the strong Szeg\"o limit theorem and we refer to \cite{greszeg,
szeg52} for more details and related results.

Today, the strong Szeg\"o limit theorem can be proved under much weaker conditions on $f$, see
e.g., \cite{goibra, hirs65, kac53}. For high dimensional extensions of the strong Szeg\"o limit
theorem, we refer to \cite{kok96} and the references therein.

In fact, the strong Szeg\"o limit theorem is a special case of a more general result proved by
Szeg\"o using the fact that $P_{n}[f]P_{n}$ is a Toeplitz form \cite{greszeg}. More specifically,
let $f$ be a real-valued integrable function such that $m\leq f(x) \leq M$. Then the eigenvalues
$\{\lambda_{k}^{(n)}\}_{k=1}^{n+1}$ of $P_{n}[f]P_{n}$ are contained in $[m, M]$. For any
continuous function $F$ defined on this interval, it was proved in \cite[Section 5.3]{greszeg} that
\begin{equation}\label{geneszego}
\lim_{n\to
\infty}\tfrac{1}{n+1}\sum_{k=1}^{n+1}F(\lambda_{k}^{(n)})=\tfrac{1}{2\pi}\int_{0}^{2\pi}F(f(x))\,
dx.
\end{equation}
Notice that~\eqref{oriszeg} is a specific case of~\eqref{geneszego}, when $F(x)=\log x$. It follows
from~\eqref{geneszego} that the bounded sequences $\{\lambda_{k}^{(n+1)}\}_{k=1}^{n+1}$ and
$\{f(\tfrac{2k\pi}{n+1})\}_{k=1}^{n+1}$ are equally distributed in the interval $[m, M]$; see
\cite[Chapter 5]{greszeg}.

By noticing that $e^{im\theta}$ is an eigenfunction of $\Delta= \tfrac{d^{2}}{dx^{2}}$, one can
view the above results as special cases of Szeg\"o limit theorem for the Laplace-Beltrami operator
(or more generally for pseudodifferential operators) on manifolds \cite{guok96, lasa96, kok96,
wid79}. \newline \noindent In the present paper, we prove the analogue of the strong Szeg\"o limit
theorem on the Sierpi\'nski gasket ($\SG$). This set is an example of fractal on which a well
established theory of Laplacian exists \cite{Ba98, Ki93, Str99, Str06}.  In this fractal setting,
the non-periodicity of the eigenfunctions of the Laplacian implies that the analogue of the matrix
$P_j[f]P_j$ is no longer related to a  Toeplitz form. Thus our results no longer follow from any
known proof of~\eqref{oriszeg}, but rather rely on the existence of localized eigenfunctions for
the Laplacian on $\SG$ \cite{BaKi97}.

The paper is organized as follows: In Section~\ref{anasg} we briefly introduce some key notions
from analysis on fractals and give a precise description of the Dirichlet spectrum of the Laplacian
on $\SG$. In Section~\ref{singleeigen} we prove a special case of the strong Szeg\"o theorem that
we use in Section~\ref{geneszeg} to prove an analogue of~\eqref{oriszeg}. Finally,
Section~\ref{lastsec} contains a further extension of the strong Szeg\"o limit theorem on $\SG$.

\section{Analysis on the Sierpi\'ski gasket}\label{anasg}
\subsection{Basic features}
In this section we collect some key facts from analysis on $\SG$ that we need to state and prove
our results.  These come from Kigami's theory of analysis on fractals, and may be found
in~\cite{Ki93}. An elementary exposition may be found in~\cite{Str99, Str06}.

Let $F_1, F_2$ and $F_3$ be the contractions defined on $\R^2$ respectively by $F_1(x) =
\frac{1}{2}x$, $F_2(x) = \frac{1}{2}x + (\frac{1}{2}, 0)$ and $F_3(x) = \frac{1}{2}x +
(\frac{1}{4}, \frac{\sqrt{3}}{4})$. The Sierpi\'nski gasket denoted by $\SG$ is the unique nonempty
compact subset of $\R^2$ such that $ \SG= \cup_{i=1}^{3}F_{i}(\SG)$.  Alternatively, $\SG$ can be
defined as a limit of graphs. For a word $\omega=(\omega_1, \omega_2, \hdots, \omega_m)$ of length
$m$, the set $F_{\omega}(\SG) =F_{\omega_{1}}\circ \cdots \circ F_{\omega_{m-1}}\circ
F_{\omega_m}(\SG)$
  with $\omega_i \in \{1, 2, 3\}$, is called an $m$-cell.
  Let $V_{0}=\{(0, 0), (1, 0),
(\frac{1}{2}, \frac{\sqrt{3}}{2})\}$, be the boundary of $\SG$ and  $V_{n}= \cup_{i=1}^{3}
F_{i}V_{n-1}$, $n\geq 1$. Define a sequence of graphs $\Gamma_{m}$ with vertices in $V_{m}$ and
edge relation $x\sim_{m} y$ given inductively by: $\Gamma_{0}$ is the complete graph with vertices
in $V_{0}$, and $x\sim_{m}y$ if and only if $x \,\mbox{and}\, y$ belong to the same $m$-cell
$F_{\omega}(\SG)$.

In all that follows, we assume that $\SG$ is equipped with the probability measure $\mu$ that
assigns the measure $3^{-m}$ to each $m-$cell.  We will also need the energy or Dirichlet form that
is naturally defined on $\SG$ and denoted $\mathcal{E}$. The precise definition of $\mathcal{E}$
will not be given here but can be found in \cite{Ki93, Str06}. All we will need in the sequel is
that $\mathcal{E}$ gives rise to a natural distance on $\SG$ called {\it the effective resistance
metric} on $\SG$, and defined for $x,y \in \SG$ by
\begin{equation}\label{resitmet}
d(x, y) = \biggparen{\min\set{\en{u}: u(x) = 0 \, \mbox{and}\, u(y) =1}}^{-1}.
\end{equation}
It is known that $d(x,y)$ is bounded above and below by constant multiples of
$|x-y|^{\log(5/3)/\log2}$, where $|x-y|$ is the Euclidean distance.

For any integer $N>1$ we will consider a partition of the Sierpi\'nski gasket ($\SG$) into
\begin{equation}\label{partition}
\SG=\cup_{|\omega|=N}F_{\omega}\SG ,\end{equation} where for each word $\omega$ of length $N$,
$F_{\omega}(\SG)$ is an $N-$cell. Note that there are exactly $3^N$ such cells, each of which has a
measure $3^{-N}$.  If $u$ is a function having support entirely contained in a single $N$-cell then
we will say that $u$ is {\em localized at scale $N$}.

\subsection{The Laplacian and its spectrum}

A Laplacian can be defined on $\SG$ either through a weak formulation using the energy and measure
or as a renormalized limit of graph Laplacians in the following manner~\cite{Ki93, Str06}. Define
the graph Laplacian $\Delta_{m}$ on $\Gamma_{m}$ by
\begin{equation}\label{graphlap}
\Delta_{m}f(x) = \sum_{y\sim_{m}x}f(y) - 4f(x)
\end{equation} for $x \in V_{m}\setminus V_0$. The Laplacian
on $\SG$ can now be defined by
\begin{equation}\label{lapla}
\Delta = \tfrac{3}{2} \lim_{m \to \infty} 5^{m} \Delta_{m}.
\end{equation}

A complete description of the spectrum of $\Delta$ on $\SG$ was given in~\cite{FuSh92} using the
method of spectral decimation introduced in~\cite{RaTo83}, and a description of the eigenfunctions
was obtained by this method in~\cite{DSV}, see also~\cite{Tep98}.   In a nutshell, the spectral
decimation method completely determines the eigenvalues and the eigenfunctions of $\Delta$ on $\SG$
from the eigenvalues and eigenfunctions of the graph Laplacians $\Delta_m$. More specifically, for
every Dirichlet eigenvalue $\lambda$ of $\Delta$ on $\SG$, there exists an integer $j\geq 1$,
called the {\em generation of birth}, such that if $u$ is a $\lambda$-eigenfunction and $k\geq j$
then $u\evald{V_{k}}$ is an eigenfunction of $\Delta_{k}$ with eigenvalue $\gamma_{k}$.  The only
possible initial values $\gamma_{j}$ are $2$, $5$ and $6$, and subsequent values can be obtained
from
\begin{equation}\label{spectdecrelnforevals}
    \gamma_{k} = \frac{1}{2}\bigl( 5+ \epsilon_{k} \sqrt{25-4\gamma_{k-1}} \bigr)\text{ for }k>j
    \end{equation}
where $\epsilon_{k}$ can take the values $\pm1$.  The sequence $\gamma_{k}$ is related to $\lambda$
by
\begin{equation}\label{spectdecrelnofgammaktolambda}
    \lambda = \tfrac{3}{2}\lim_{k\to \infty} 5^{k}\gamma_{k}.
    \end{equation}
In particular the convergence of~\eqref{spectdecrelnofgammaktolambda} implies that $\epsilon_{k}=1$
for at most a  finite number of $k$ values.  We let $l=\min\{k:\epsilon_{k}=-1\}$ and call it the
{\em generation of fixation}.

An interesting and useful feature is that there are a great many eigenfunctions which satisfy both
Dirichlet and Neumann boundary conditions; this is a general property of the Laplacian on fractals
with sufficient symmetry~\cite{BaKi97}, and in the case of $\SG$ it implies both that most
eigenfunctions are localized on small sets, and that eigenspaces have high multiplicity.

Using the spectral decimation algorithm and elementary properties of the map
in~\eqref{spectdecrelnforevals} one can see that the size of an eigenvalue depends (up to constant
factors) on its generation of fixation, and its multiplicity depends on its generation of birth. We
summarize the relevant features of eigenvalues at the bottom of the spectrum and their eigenspaces
in Proposition~\ref{evalsefnsetc} below, and refer the reader to the original treatments
~\cite{BaKi97, DSV, FuSh92}, or the exposition in~\cite{Str06} for proofs.

\begin{prop}\label{evalsefnsetc}
There is a constant $\kappa$ such that the $\tfrac{1}{2}(3^{m+1}-3)$ smallest eigenvalues of
$-\Delta$ are precisely those with size at most $\kappa5^{m}$, and all have generation of fixation
$l\leq m$.  The eigenvalues, their multiplicities, and bases for their eigenspaces may be described
as follows.
\begin{itemize}
    \item The $2$-series eigenvalues are those obtained from~\eqref{spectdecrelnforevals}
and~\eqref{spectdecrelnofgammaktolambda} with generation of birth $j=1$ and $\gamma_{j}=2$.  Each
such eigenvalue has multiplicity $1$.
    \item The $5$-series eigenvalues are those obtained from~\eqref{spectdecrelnforevals}
and~\eqref{spectdecrelnofgammaktolambda} with any generation of birth $j\geq1$ and $\gamma_{j}=5$.
There are $2^{m-j}$ such eigenvalues for each $1\leq j\leq m$, every one having multiplicity
$\tfrac{1}{2}(3^{j-1}+3)$.   For each such eigenvalue, there is a basis for the corresponding
eigenspace in which all but two of the basis functions have support in a collection of
$(j-1)$-cells arranged in a loop around a ``hole'' of scale at least $(j-1)$ in $\SG$.  There are
$\tfrac{1}{2}(3^{j-1}-1)$ such holes and one eigenfunction for each hole.
    \item The $6$-series eigenvalues are those obtained from~\eqref{spectdecrelnforevals}
and~\eqref{spectdecrelnofgammaktolambda} with any generation of birth $j\geq2$, $\gamma_{j}=6$ and
$\epsilon_{j+1}=+1$.  There are $2^{m-j-1}$ such eigenfunctions for each $2\leq j<m$ and $1$ for
$j=m$, every one having multiplicity $\tfrac{1}{2}(3^{j}-3)$.  For each such eigenvalue there is a
basis for the corresponding eigenspace that is indexed by points of $V_{j-1}\setminus V_{0}$, and
in which each basis element is supported on the union of the two $j$-cells that intersect at the
corresponding point in $V_{j-1}\setminus V_{0}$.
\end{itemize}
\end{prop}

Recall that a function is said to be localized at scale $N$ if its support is contained entirely in
a single $N$-cell.  It is apparent from the above that there are eigenfunctions that are localized
at scale $N$ provided $j>N$.  For later use, we compute the number of these that occur in each of
the $5$ and $6$-series.

Let $j>N$ and consider a $6$-series eigenvalue with generation of birth $j$. In the associated
eigenspace there are $\tfrac{1}{2}(3^{N+1}-3)$ basis functions corresponding to the vertices in
$V_{N}\setminus V_{0}$, and which are not localized at scale $N$. The remaining
$\tfrac{1}{2}(3^{j}-3^{N+1})$ basis functions correspond to vertices in $V_{j-1}\setminus V_{N}$
and are localized at scale $N$.

For $j>N$ and a fixed $5$-series eigenvalue, the basis elements for the eigenspace are supported on
loops or chains. A loop is contained in a cell of scale $N$ if and only if the hole it encircles
has scale at least $N+1$. The number of holes of scale at most $N$ is $\tfrac{1}{2}(3^{N}-1)$, so
the number of basis elements not localized at scale $N$ is this plus the two that are not loops,
giving $\tfrac{3^{N}}{2}$ in total.  The remaining $\tfrac{1}{2}(3^{j-1}-3^{N})$ basis
eigenfunctions are localized at scale $N$.

\section{Szeg\"o limit theorem on $\SG$ for a single eigenspace}\label{singleeigen}

In this section we prove a Szeg\"o limit theorem for a single $5$-series or $6$-series eigenspace
of the Laplacian on $\SG$.

Let $\lambda_{j}$ be a $6$-series eigenfunction with generation of birth $j$ and eigenspace
$E_{j}$. Denote the span of those eigenfunctions corresponding to $\lambda_j$ that are localized at
scale $N<j$ by $E_{j}^{N}$. Let $d_{j}^{N}=\dim{E_{j}^{N}} = \tfrac{1}{2}(3^{j} -3^{N+1})$ and
$\alpha_{j}^{N}=\tfrac{1}{2}(3^{N+1}-3)$ be the dimension of the complementary space in $E_{j}$.
Since there are $3^{N}$ cells of scale $N$, we see that the number of eigenfunctions supported on a
single cell is $m_{j}^{N}=\tfrac{1}{2}(3^{j-N}-3)$.

An analogous construction may be done for a $5$-series eigenfunction, with the only change being
that in this case $\alpha_{j}^{N}=\tfrac{3^{N}-3}{2}$, $d_{j}^{N}=\tfrac{1}{2}(3^{j-1}-3^{N})$ and
$m_{j}^{N}=\tfrac{1}{2}(3^{j-N-1}-1)$.

For each $N$-cell, use the Gram-Schmidt process to orthonormalize the collection of eigenfunctions
supported on that cell.  Since functions on separate cells are already orthogonal, taking the union
over all $N$-cells gives an orthonormal basis $\{\tilde{u}_{k}\}_{k=1}^{d_{j}^{N}}$ for
$E_{j}^{N}$. Adjoining the remaining basis elements of $E_{j}$ and again using the Gram-Schmidt
process extends this to an orthonormal basis
\begin{equation*}
\{u_{k}\}_{k=1}^{d_{j}}=\{\tilde{u}_{k}\}_{k=1}^{d_{j}^{N}}\cup \{v_{k}\}_{k=1}^{\alpha_{j}^{N}},
\end{equation*}
for $E_{j}$, where only the $v_{k}$ are not localized at scale $N$.

Let $P_{j}$ be the projection of $L^{2}(\SG)$ onto $E_{j}$.  For $g \in L^{2}(\SG)$, $P_j$ is
defined by
\begin{equation}\label{defproj}
P_{j}g(x):=\sum_{k=1}^{d_{j}}g_{k}u_{k}(x)=\sum_{k=1}^{d_{j}}\ip{g}{u_{k}}u_{k}(x).
\end{equation}
For a real-valued function $f$ on $\SG$, we recall that $[f]$ the operator corresponding to the
pointwise multiplication by $f$.

\subsection{The case of simple functions}\label{subsec1}

\begin{lemma}\label{simpleszeg}
Let $f=\sum_{k=1}^{3^{N}}a_{k}\chi_{C_{k}}$ where we assume that $a_{k}>0$ for  all $k$ and
$\chi_{C_{k}}$ denotes the characteristic function of the $N$-cell $C_{k}$. Then for $P_{j}$ as
above,
$$\lim_{j\to \infty}\tfrac{1}{d_{j}}\log \det P_{j}[f]P_{j}=\int_{\SG}\log f(x)\, d\mu(x).$$
Furthermore, for $j$ large enough,
$$\tfrac{1}{d_{j}}\log \det P_{j}[f]P_{j} - \int_{\SG} \log f(x)\, d\mu(x) = O(d_{j}^{-1}).$$
\end{lemma}

\begin{proof}
For $j>N$ set $M_{j}=P_{j}[f]P_{j}$.  Then $M_{j}$ is a $d_{j}\times d_{j}$ matrix which has block
structure
\begin{equation}\label{matrixmj}
    M_{j}=\begin{bmatrix} R_{j}& \star \\
    0 & N_{j}\end{bmatrix}
\end{equation}
with respect to the basis $\{u_{k}\}_{k=1}^{d_{j}}$. Here $R_j$ is an invertible $d_{j}^{N}\times d^{N}_{j}$ matrix
corresponding to the ``localized'' part, while $N_{j}$ is an invertible $\alpha_{j}^{N} \times
\alpha_{j}^{N}$ matrix corresponding to the ``non-localized'' part. Furthermore, $R_j$ is a block
diagonal matrix, where each block is an $m_{j}^{N}\times m_{j}^{N}$ matrix  that corresponds to a
single $N-$cell $C_{k}$, and is therefore simply $a_{k}I_{m_{j}^{N}}$ where $I_{d}$ is notation for
the $d\times d$ identity matrix. It follows immediately that
\begin{align*}
    \log\det M_{j}
    &= \log\det R_{j} +\log\det N_{j}\\
    &= \log \biggl(\prod_{k=1}^{3^{N}} a_{k}^{m_{j}^{N}}\biggr)  +\log\det N_{j}\\
    &= m_{j}^{N} \biggl( \sum_{k=1}^{3^{N}} \log a_{k} \biggr) +\log\det N_{j}.
    \end{align*}
Using $d_{j}^{N}=3^{N}m_{j}^{N}$ we conclude
\begin{align*}
    \tfrac{1}{d_{j}}\log \det M_{j}
    &=\tfrac{m_{j}^{N}3^{N}}{d_{j}} \sum_{k=1}^{3^{N}}3^{-N}\log a_{k} +  \tfrac{1}{d_{j}}\log \det N_{j} \notag \\
    &=\tfrac{d^{N}_{j}}{d_{j}} \int_{\SG}\log f (x)\, d\mu(x) +  \tfrac{1}{d_{j}}\log \det N_{j},
    \end{align*}
and since $d_{j}-d_{j}^{N}=\alpha_{j}^{N}$, we have
\begin{align}\label{errorterminsimplecase}
    \lefteqn{\tfrac{1}{d_{j}}\log \det M_{j} -  \int_{\SG}\log f (x)\, d\mu(x)}\quad&\notag\\
    &= \tfrac{-\alpha_{j}^{N}}{d_{j}}\int_{\SG}\log f (x)\, d\mu(x) + \tfrac{1}{d_{j}}\log\det N_{j}.
    \end{align}

We can now afford a crude estimate of the term $\log\det N_{j}$. Since $f\in L^{\infty}$, the
multiplier $[f]$ is bounded on $L^{2}$ by $\|f\|_{\infty}$.  It follows that $\langle
N_{j}g,g\rangle \leq \|f\|_{\infty}$ for any $g\in E_{j}$ with $\|g\|_{2}=1$, and therefore that
$\det N_{j}\leq \|f\|_{\infty}^{\alpha_{j}^{N}}$.  Combining this
with~\eqref{errorterminsimplecase} we see
\begin{equation*}
    \biggl| \tfrac{1}{d_{j}}\log \det M_{j} -  \int_{\SG}\log f (x)\, d\mu(x) \biggr|
    \leq \tfrac{\alpha_{j}^{N}}{d_{j}} \Bigl( \|\log f (x)\|_{1} + \|f\|_{\infty} \Bigr)
    \end{equation*}
which completes the proof because $\alpha_{j}^{N}$ is bounded by a constant multiple of $3^{N}$ and
$d_{j}$ is comparable to $3^{j}$.
\end{proof}

\subsection{The case of continuous positive functions}\label{contszeg}
\begin{thm}\label{main1}
 Let $f$ be a positive and continuous function on $\SG$. Then
 \begin{equation}\label{mainest1}
 \lim_{j\to \infty}\tfrac{1}{d_{j}}\log \det P_{j}[f]P_{j}=\int_{\SG}\log f(x)\, d\mu(x).
 \end{equation}
If in addition, we assume that $f$ is H\"older continuous of order $\alpha$ in the resistance
metric $R$ on $\SG$, then
\begin{equation}\label{improvedest1}
\tfrac{1}{d_{j}}\log \det P_{j}[f]P_{j} - \int_{\SG}\log f(x)\, d\mu(x)  = O(d_{j}^{-\beta})
\end{equation} where
$$\beta =\tfrac{\alpha \log(5/3)}{\log 3 + \alpha \log(5/3)}=1 - \tfrac{\log 3}{\alpha \log(5/3) + \log 3}.$$
\end{thm}

\begin{proof}
Since $\SG$ is compact, $\min_{x\in\SG} f(x)=m>0$. Given $\epsilon>0$, uniform continuity provides
$N$ and a simple function $f_{N}=\sum_{k=1}^{3^{N}}a_{k}\chi_{C_{k}}$ such that
\begin{equation*}
    \|f-f_{N}\|_{\infty}<\min\bigl(\tfrac{1}{2},\tfrac{m}{2},\tfrac{\epsilon m}{2}\bigr),
    \end{equation*}
from which the following are immediate,
\begin{gather}
    \tfrac{\bigl| f(x)-f_{N}(x)\bigr|}{\bigl|f_{N}(x)\bigr|}
    \leq \epsilon, \notag\\
    1-\epsilon \leq \tfrac{f(x)}{f_{N}(x)} \leq 1+\epsilon,\label{ratioboundsinmain1thm}\\
    -2\epsilon\leq \log(1-\epsilon) \leq \log \Bigl( \tfrac{f(x)}{f_{N}(x)}\Bigr) \leq
    \log(1+\epsilon)\leq \epsilon. \label{logboundsinmain1thm}
    \end{gather}
Note that~\eqref{logboundsinmain1thm} implies $\bigl| \int_{SG}\log f-\int_{SG}\log f_{N}\bigr|\leq
2\epsilon$.

Now let us estimate $\log\det P_{j}[f]P_{j}$ in the same manner as was done in
Lemma~\ref{simpleszeg}.  It has a block structure like~\eqref{matrixmj}, but the diagonal blocks in
$R_{j}$ are no longer multiples of the identity matrix $I_{m_{j}^{N}}$.  However it follows
from~\eqref{ratioboundsinmain1thm} that the values on the diagonal corresponding to $C_{k}$ are
bounded below by $a_{k}(1-\epsilon)$ and above by $a_{k}(1+\epsilon)$, and thus
\begin{equation*}
    \tfrac{d_{j}^{N}}{d_{j}} \log(1-\epsilon)
    \leq \tfrac{1}{d_{j}} \log\det R_{j} -  \tfrac{d_{j}^{N}}{d_{j}} \int_{\SG} f_{N} \, d\mu
    \leq \tfrac{d_{j}^{N}}{d_{j}} \log(1+\epsilon)
    \end{equation*}
and in particular $\bigl| d_{j}^{-1}\log\det R_{j} - d_{j}^{N}d_{j}^{-1} \int_{SG} f_{N} \bigr|\leq
2\epsilon$.

Combining these estimates with the same $\log\det N_{j} \leq \alpha_{j}^{N}\|f\|_{\infty}$ bound
used in Lemma~\ref{simpleszeg} we have
\begin{align}
    \lefteqn{\biggl| \tfrac{1}{d_{j}} \log\det P_{j}[f]P_{j} - \int_{SG} \log f\, d\mu\biggr| }\quad& \notag\\
    &\leq \biggl| \tfrac{1}{d_{j}} \log\det R_{j} - \tfrac{d_{j}^{N}}{d_{j}} \int_{SG} \log f_{N} \, d\mu \biggr|
        + \tfrac{\alpha_{j}^{N}}{d_{j}} \biggl| \int_{SG} \log f_{N} \, d\mu \biggr| \notag\\
    &\quad + \biggl| \int_{SG} \log f_{N}- \log f \, d\mu \biggr| + \tfrac{1}{d_{j}}\log\det N_{j} \notag\\
    &\leq 4\epsilon + c\tfrac{\alpha_{j}^{N}}{d_{j}}\bigl( \|f\|_{\infty}+\|\log f\|_{1} \bigr) \notag\\
    &\leq 4\epsilon + c3^{N-j}\bigl( \|f\|_{\infty}+\|\log f\|_{1} \bigr). \label{mainthm1boundforall}
    \end{align}
This gives the first statement of the theorem.

In the case that $f$ is H\"older continuous of order $\alpha$ in the resistance metric we see that
$$\epsilon = \nm{f-f_{N}}_{L^{\infty}} = O((\tfrac{3}{5})^{N\alpha}).$$  For a fixed large
$j$ we may then choose $N$ such that the bound in~\eqref{mainthm1boundforall} is minimized, which
occurs when $\epsilon\approx3^{N-j}.$  Setting $(3/5)^{N\alpha}=3^{N-j}$ we compute
\begin{equation*}
    3^{N-j} =3^{\tfrac{-j\alpha \log(5/3)}{\log 3 + \alpha \log(5/3)}}
    \end{equation*}
and substitute into~\eqref{mainthm1boundforall} to obtain~\eqref{improvedest1}, using
$d_{j}\approx3^{j}$.
\end{proof}
\begin{rem}
The following special cases of~\eqref{improvedest1} are worth pointing out. \newline If $\alpha =
1/2$ (which is the case when $f\in dom \mathcal{E}$), then $\beta = 1 -\tfrac{\log 9}{\log 15}.$
\newline If $\alpha =1$ (which is the case if $f \in dom \Delta$), then $\beta = 1-\tfrac{\log
3}{\log 5}.$
\end{rem}

\section{General Szeg\"o Theorem on $\SG$}\label{geneszeg}
In this section we prove analogues of the results proved in Section~\ref{singleeigen} for the
situation where we look at all eigenvalues up to a certain value $\Lambda$.  We therefore let
$E_{\Lambda}$ be the span of all eigenfunctions corresponding to eigenvalues $\lambda$ of $-\Delta$
for which $\lambda\leq \Lambda$, let $P_{\Lambda}$ be projection onto $E_{\Lambda}$, and set
$d_{\Lambda}=\dim(E_{\Lambda})$. We also suppose that a scale $N$ is fixed.

Since $-\Delta$ is self-adjoint the eigenspaces of distinct eigenvalues are orthogonal.  For each
$\lambda<\Lambda$ from either the $5$-series or the $6$-series, and having generation of birth
$j>N$ we take an orthonormal basis for the corresponding eigenspace of the type described in
Section~\ref{singleeigen}.  For all other eigenspaces in $E_{\Lambda}$ we simply take orthonormal
bases.  The union of the basis vectors is then a basis for $E_{\Lambda}$ and in this basis the
operator $P_{\Lambda}[f]P_{\Lambda}$ is a block diagonal matrix $M_{\Lambda}$ with one block
$M_{\lambda}$ for each eigenvalue $\lambda\leq\Lambda$.

\begin{thm}\label{main2} Let $f> 0$ be a continuous function on $\SG$. Then,
\begin{equation}\label{mainest2}
\lim_{\Lambda\to \infty}\tfrac{1}{d_{\Lambda}}\log \det M_{\Lambda} =\int_{\SG}\log f(x)\, d\mu(x).
\end{equation}
If in addition we assume that $f$ is H\"older continuous of order $\alpha$ in the resistance metric
$R$ on $\SG$ then
\begin{equation}\label{improvedest2}
    \tfrac{1}{d_{\Lambda}}\log \det M_{\Lambda}-\int_{\SG}\log f(x)\, d\mu(x)
    =O(d_{\Lambda}^{-\tilde{\beta}})
\end{equation}
where $$\tilde{\beta} =\beta (1-\tfrac{\log 2}{\log 3})=(\tfrac{\alpha \log(5/3)}{\log 3 + \alpha
\log(5/3)})(1-\tfrac{\log 2}{\log 3}).$$
\end{thm}

\begin{proof}
Fix $\epsilon>0$.  It is clear that $\log\det M_{\Lambda} = \sum_{\lambda\leq\Lambda}\log\det
M_{\lambda}$. If $\lambda$ is one of the $5$ or $6$-series eigenvalues with $j>N$ then replacing
$f$ by $f_{N_{0}}$ as in the proof of Theorem~\ref{main1} we have from~\eqref{mainthm1boundforall}
\begin{equation}\label{estimateforgoodeigenspaceblocks}
    \biggl| \log\det M_{\lambda} - d_{\lambda}\int_{SG} \log f\, d\mu \biggr|
    \leq 4\epsilon d_{\lambda} + c3^{N}\bigl( \|f\|_{\infty}+\|\log f\|_{1} \bigr).
    \end{equation}
Now we let $\Gamma_{N}$ be the set of $\lambda<\Lambda$ with generation of birth $j>N$ and
sum~\eqref{estimateforgoodeigenspaceblocks} over $\lambda\in\Gamma_{N}$, noting that
$d_{\Lambda}^{-1}\sum_{\lambda\in\Gamma_{N}}d_{\lambda}\leq 1$. Using the trivial bound $\log\det
M_{\lambda}\leq d_{\lambda} \|f\|_{\infty}$ for the remaining terms in $\log\det M_{\Lambda}$, and
$d_{\lambda}\|\log f\|_{1}$ for those making up the integral, we obtain
\begin{align}
    \lefteqn{\biggl| \tfrac{1}{d_{\Lambda}} \log\det M_{\Lambda} - \int_{SG} \log f\, d\mu
    \biggr|}\quad&\notag\\
    &\leq 4\epsilon + \bigl( \|f\|_{\infty}+\|\log f\|_{1} \bigr)
    \biggl( \Bigr( \tfrac{c3^{N}}{d_{\Lambda}} \sum_{\lambda\in\Gamma_{N}} 1 \Bigr)+
    \Bigl(\sum_{\lambda\not\in\Gamma_{N}} \tfrac{d_{\lambda}}{d_{\Lambda}} \Bigr) \biggr).
    \label{main2intermedest}
    \end{align}

The remaining work in the proof is to estimate the number of eigenvalues in $\Gamma_{N}$ and the
sum of the dimensions $d_{\lambda}$ for $\lambda\not\in\Gamma_{N}$, which we do using
Proposition~\ref{evalsefnsetc}. For this purpose, take $m\in\mathbb{N}$ so
$\kappa5^{m-1}\leq\Lambda<\kappa5^{m}$, where $\kappa$ is as in Proposition~\ref{evalsefnsetc}.
Since $\Gamma_{N}$ is empty and the estimate is trivial if $N\geq m$, we assume without loss of
generality that $N<m$. The eigenvalues less than $\kappa5^{m}$ and having generation of birth
$j\leq N$ number $2^{m-1}$ with multiplicity $1$ from the $2$-series, $2^{m-j}$ with multiplicity
$\tfrac{1}{2}(3^{j-1}+3)$ from the $5$-series and $2^{m-j-1}$ with multiplicity
$\tfrac{1}{2}(3^{j}-3)$ if $2\leq j\leq N$. Summing these gives
\begin{equation*}
    2^{m-1} + \sum_{1}^{N} 2^{m-j-1}(3^{j-1}+3) + \sum_{2}^{N} 2^{m-j-2}(3^{j}-3)
    =O(2^{m-N}3^{N}).
    \end{equation*}
so
\begin{equation}\label{main2intermedest2}
    \sum_{\{\lambda\not\in\Gamma_{N}: \lambda\leq \Lambda\}} d_{\lambda}
    \leq \sum_{\{\lambda\not\in\Gamma_{N}: \lambda\leq \kappa 5^{m+1}\}} d_{\lambda}
    =O(2^{m-N}3^{N}).
    \end{equation}
Moreover the number of $\lambda\in\Gamma_{N}$ such that $\lambda\leq\kappa 5^{m}$ is
$\sum_{N+1}^{m}2^{m-j}$ from the $5$ series and $\sum_{N+1}^{m} 2^{m-j-1}$ from the $6$ series,
giving a total that is $O(2^{m-N})$.  This implies
\begin{equation*}
    \tfrac{c3^{N}}{d_{\Lambda}} \sum_{\lambda\in\Gamma_{N}} 1
    =O( 2^{m-N} 3^{N}d_{\Lambda}^{-1}),
    \end{equation*}
and substituting this and~\eqref{main2intermedest2} into~\eqref{main2intermedest}, along with
$d_{\Lambda}\geq \tfrac{1}{2}(3^{m}-3)$ because $\Lambda\geq\kappa5^{m-1}$, we have
\begin{equation*}
    \biggl| \tfrac{1}{d_{\Lambda}} \log\det M_{\Lambda} - \int_{SG} \log f\, d\mu \biggr|
    \leq 4\epsilon + c\Bigl(\tfrac{3}{2} \Bigr)^{N-m},
    \end{equation*}
which proves the first statement of the theorem.

For H\"{o}lder continuous $f$ and fixed $\Lambda$ we may now optimize the choice of $N$ as in the
proof of Theorem~\ref{main1} to obtain~\ref{improvedest2}.
\end{proof}

\begin{rem}
Observe that in comparison with~\eqref{improvedest1}, the error in~\eqref{improvedest2} is decaying
at a slower rate. This is a consequence of the fact that, at the optimal $N$, the eigenfunctions
that are not localized at scale $N$ make up a larger proportion (in terms of dimension) of the
space $E_{\Lambda}$ than they do in the spaces $E_{\lambda}$ with $\lambda\approx\Lambda$.
\end{rem}

\section{``Almost'' equally distributed sequences}\label{lastsec}

As mentioned to in the Introduction,~\eqref{geneszego} can be translated into results on equally
distributed sequences. In this section we shall prove an analogue of~\eqref{geneszego} on $\SG$.
This will be used to define the notion of ``almost'' equally distributed sequences on $\SG$. We
recall the definition of equally distributed sequences due to H.~Weyl, for which we refer to
\cite[Chapter 5]{greszeg}.
\begin{deft}\label{eqseq}
Fix $K>0$.  For each $n$ let $a_{1}^{(n)},\dotsc, a_{n+1}^{(n)}$ and $b_{1}^{(n)},\dotsc,
b_{n+1}^{(n)}$ be sets of $n+1$ numbers from the interval $[-K,K]$.  We say that the sets
$\{a_{j}^{(n)}\}$ and $\{b_{k}^{(n)}\}$, $n\to \infty$, are equally distributed in the interval
$[-K, K]$, if given any continuous function $F$ on $[-K, K]$ we have
\begin{equation*}
    \lim_{n \to \infty} \tfrac{1}{n+1} \sum_{j=1}^{n+1}[F(a_{j}^{(n)}) - F(b_{j}^{(n)})]=0.
    \end{equation*}
\end{deft}

We first consider an extension of the results proved in Section~\ref{singleeigen}.  Recall that
$M_j = P_{j}[f]P_{j}$, where $P_j$ is the orthogonal projection onto the eigenspace $E_{j}$
corresponding to a $5$ or $6$-series eigenvalue $\lambda_{j}$ of $-\Delta$ with generation of birth
$j$. Let $\bigl\{\sigma_{k}^{(j)}\bigr\}_{k=1}^{d_{j}}$ be the eigenvalues of $M_{j}$

\begin{lemma}\label{eqseq1}
Let $f=\sum_{k=1}^{3^{N}}a_{k}\chi_{C_{k}}$ with all $a_{k}>0$ and let $m=\min_{k}a_{k}$,
$M=\max_{k}a_{k}$. Let $F$ be continuous on $[m,M]$.  Then
\begin{equation*}
    \lim_{j\to \infty}\tfrac{1}{d_{j}}\sum_{k=1}^{d_{j}}F(\sigma_{k}^{(j)})=\int_{\SG}F(f(x))\,
    d\mu(x).
\end{equation*}
Moreover there is a set of points $\{s^{(j)}_{k}\}_{k}^{d_{j}}$ in $\SG$ such that
$\{\sigma_{k}^{(j)}\}$ and $\{f(s_{k}^{(j)})\}$ are ``almost'' equally distributed in $[m, M]$ when
$j \to \infty$.
\end{lemma}

\begin{rem}
We use the term ``almost'' equally distributed because the above limit is computed along the
subsequence $d_{j}$ of the positive integers.
\end{rem}

\begin{proof}
Note from the proof of Lemma~\ref{main1} that all eigenvalues of $M_{j}$ satisfy $m \leq
\sigma_{k}^{(j)} \leq M$, and that each $a_k$ is an eigenvalue of $M_j$ with multiplicity
$m_{j}^{N}$. The remaining eigenvalues form a set $\Upsilon_{j}$ with $\#\Upsilon_{j}\leq
\alpha_{j}^{N}$.  Using the fact that $m_{j}^{N}=3^{-N}d_{j}^{N}$ we compute
\begin{align*}
    \lefteqn{\biggl|\tfrac{1}{d_{j}} \sum_{k=1}^{d_{j}}F(\sigma_{k}^{(j)}) - \int_{SG} F(f(x))\, d\mu
    \biggr|}\quad&\\
    &= \biggl| \tfrac{1}{d_{j}} \sum_{\Upsilon_{j}} F(\sigma_{k}^{(j)})
        + \tfrac{m_{j}^{N}}{d_{j}} \sum_{k=1}^{3^{N}}F(a_{k})  - \int_{SG} F(f(x))\, d\mu \biggr|\\
    &\leq \tfrac{\alpha_{j}^{N}}{d_{j}}\|F\|_{L^{\infty}([m,M])} + \Bigl(\tfrac{d_{j}^{N}}{d_{j}}-1\Bigr) \int_{SG} F(f(x))\, d\mu \\
    &\leq \tfrac{2\alpha_{j}^{N}}{d_{j}}\|F\|_{L^{\infty}([m,M])}.
    \end{align*}
The last part of the result follows from approximating $\int_{\SG}F(f(x))\, d\mu(x)$ with Riemann
sums.
\end{proof}

With this result we can prove the following extension of Theorem~\ref{main1}

\begin{thm}\label{maineqseq1}
Let $f>0$ be continuous on $\SG$ and $M=\max_{\SG}f(x)$. If $F$ is continuous on $[0, \infty)$,
then
$$\lim_{j\to \infty}\tfrac{1}{d_{j}}\sum_{k=1}^{d_{j}}F(\sigma_{k}^{(j)})=\int_{\SG}F(f(x))\,
d\mu(x).$$  Moreover there is a set of points $\{s^{(j)}_{k}\}_{k}^{d_{j}}$ in $\SG$ such that
$\{\sigma_{k}^{(j)}\}$ and $\{f(s_{k}^{(j)})\}$ are ``almost'' equally distributed in $[0,
\nm{f}_{L^{\infty}}]$ as $j\to\infty$.
\end{thm}

\begin{proof}
Let $\delta>0$ be given and let $0<\epsilon<\delta$ be such that $|a-b|\leq\epsilon$ and
$a,b\in[0,M]$ implies $|F(a)-F(b)|<\delta$.  Take $f_{N}$ a simple function as in the proof of
Theorem~\ref{main1}.  We saw in that proof that the eigenvalues of $P_{j}[f]P_{j}$ that correspond
to the eigenfunctions localized at scale $N$ are bounded below by $(1-\epsilon)a_{k}$ and above by
$(1+\epsilon)a_{k}$.  Writing $\rho_{k}^{j}$ for the eigenvalues of $P_{j}[f_{N}]P_{j}$ we have
\begin{align*}
    \lefteqn{\biggl| \tfrac{1}{d_{j}} \sum_{k=1}^{d_{j}}F(\sigma_{k}^{(j)}) - \int_{SG} F(f(x))\, d\mu
    \biggr|}\quad&\notag\\
    &\leq \biggl| \tfrac{1}{d_{j}} \sum_{k=1}^{d_{j}}F(\rho_{k}^{(j)}) - \int_{SG} F(f_{N}(x))\, d\mu \biggr|
        + \biggl| \tfrac{1}{d_{j}} \sum_{k=1}^{d_{j}}\Bigl( F(\sigma_{k}^{(j)}) - F(\rho_{k}^{(j)}) \Bigr) \biggr| \notag\\
    &\quad + \biggl| \int_{SG} \Bigl( F(f(x))- F(f_{N}(x))\Bigr) \, d\mu(x) \biggr| \notag\\
    &\leq \tfrac{2\alpha_{j}^{N}}{d_{j}}\|F\|_{L^{\infty}([m,M])} + \delta + \epsilon.
    \end{align*}
This proves the first statement, and the second statement follows using Riemann sums as before.
\end{proof}

More generally we have the following extension of Theorem~\ref{main2}, in which we denote the
eigenvalues of $M_{\Lambda}$ by $\{\sigma_{k}^{(\Lambda)}\}_{k=1}^{d_{\Lambda}}$.
\begin{thm}\label{maineqseq2}
Let $f> 0$ be continuous on $\SG$ and $F$ be continuous on $[0, \infty)$. Then
\begin{equation}\label{maineqseq2eqn}
    \lim_{\Lambda\to\infty}\tfrac{1}{d_{\Lambda}} \sum_{k=1}^{d_{\Lambda}} F(\sigma_{k}^{(\Lambda)})
    =\int_{\SG}F(f(x))\, d\mu(x).
\end{equation}
Moreover there is a set of points $\{s^{(\Lambda)}_{k}\}_{k}^{d_{\Lambda}}$ in $\SG$ such that
$\{\sigma_{k}^{(\Lambda)}\}$ and $\{f(s_{k}^{(\Lambda)})\}$ are ``almost'' equally distributed in $[0,
\nm{f}_{L^{\infty}}]$ as $\Lambda \to \infty$.
\end{thm}

\begin{proof}
For $\delta>0$ and $N$ as in Theorem~\ref{maineqseq1} we decompose the sum into terms corresponding
to $\lambda\in\Gamma_{N}$ and $\lambda\not\in\Gamma_{N}$ as in Theorem~\ref{main2}. For
$\lambda\in\Gamma_{N}$ we estimate as in the proof of Theorem~\ref{maineqseq1}, and follow the
argument of Theorem~\ref{main2} to find that the left and right sides of~\eqref{maineqseq2eqn}
differ by at most $2\delta+O(2^{m-N}3^N)$ where $\Lambda\approx5^{m}$.
\end{proof}

\begin{rem}
The results proved here for $\SG$ should extend to other fractals on which localized eigenfunctions
come to predominate in the spectrum as the eigenvalues increase.  Sufficient symmetry conditions
for the existence of high multiplicity eigenspaces with localized eigenfunctions were given
in~\cite{BaKi97}.
\end{rem}

\subsection*{Acknowledgment}  The authors are grateful to Victor Guillemin for suggesting that we
investigate these questions.

\end{document}